\documentclass[preprint,12pt]{elsarticle}

\usepackage{graphicx}
\usepackage{amssymb}
\usepackage{lineno}
\usepackage{graphicx}
\usepackage{color}
\usepackage{subfigure}
\usepackage{url}
\usepackage{amsmath,amsfonts,amssymb}
\usepackage{array,multirow}
\usepackage{graphics}
\usepackage{float}
\usepackage{caption}
\usepackage{verbatim}
\usepackage{color}
\usepackage{rotating}
\usepackage{tabularx}
\usepackage{lscape}
\usepackage[table]{xcolor}
\usepackage{tikz}
\usepackage{pgfplots}
\usepackage{pgfplotstable}
\usepackage[normalem]{ulem}
\usetikzlibrary{arrows,%
                petri,%
                topaths,
 shapes,trees}%

\usepackage{algorithmic}

\usepackage{booktabs}
\begin{document}

\begin{frontmatter}


\title{Models and algorithms for the Flying Sidekick Traveling Salesman Problem}



\author[label1]{Mauro Dell'Amico} 
%\ead{mauro.dellamico@unimore.it}
\author[label1]{Roberto Montemanni} 
%\ead{roberto.montemanni@unimore.it}
\author[label1]{Stefano Novellani} 
\ead{stefano.novellani@unimore.it}
\address[label1]{Dipartimento di Scienze e Metodi dell'Ingegneria (DISMI), \\Universit\`a di Modena e Reggio Emilia (UNIMORE)\\
            via Amendola 2, 42122 Reggio Emilia, Italy}

\begin{abstract}
This paper presents a set of new formulations for the Flying Sidekick Traveling Salesman Problem, where a truck and a drone cooperate to deliver parcels to customers minimizing the completion time.
The new formulations improve the results of the literature by solving to optimality several benchmark instances for which an optimal solution was previously unknown. A matheuristic algorithm, strongly based on the new models, is also discussed. Experimental results show that this method is able to provide good quality solutions in short time even for the larger instances, on which the mathematical models struggle to provide either good heuristic solution or strong lower bounds.

\end{abstract}


\end{frontmatter}

%\linenumbers

\section{Introduction}\label{sec:intro}
The use of aerial drones or unmanned aerial vehicles is
gaining more and more relevance in several non military fields, from precision agriculture, to logistics operations, to catastrophic events management, etc. Their use to a large and diverse set of applications is due to the advantages that result from their flexibility, agility, and usability mainly due to the small size and the fact that no human is needed on board.
Drone applications can be divided, roughly, into two: the collection of data and information, and the movimentation of goods. This last case is the one that interests us, in particular, we consider a problem linked to parcels delivery. The boom of e-commerce and the promise of faster and faster deliveries has provided a challenging task to e-commerce and express delivery companies. Among the first companies exploring the use of drones in parcel deliveries one can mention Alibaba, Alphabet, Amazon 
 and JD.com. We address the reader to a recent survey on optimization approaches for drones in the civil sector by Otto et. al. \cite{ottooptimization}.

In this paper we consider the {Flying Sidekick Traveling Salesman Problem}, a problem where parcels are delivered to customers either by a truck or (if possible) by a drone. The two vehicles work coupled such that the drone can leave and must return to the truck after visiting one customer, performing flights not exceeding its battery endurance. The synchronization among the two vehicles is essential and the completion time at the end of the operations should be minimized.
We propose some enhanced mixed integer linear programming (MILP) formulations for the problem, in particular: a three indexed and a pair of two indexed formulations. For two of the formulations we introduce a set of variables that account for the drone presence on the truck -- this helps in avoiding infeasible drone flights -- while for the third one we consider this explicitly with a set of constraints, the element of which are separated in a branch-and-cut fashion. Another novelty with respect to previous works is that in all the proposed formulation we use only one set of timing variables for both the truck and the drone to account for synchronization. We also introduce a simple matheuristic method strictly based on the new formulations, which is shown work well especially on larger instances, for which the formulations encounter scalability issues.

{
The paper is organized as follows: in Section \ref{sec:lit} we discuss upon the most relevant related literature, while in Section \ref{sec:prob} we formally describe the problem. New mathematical formulations and their implementations are described in Sections \ref{sec:form}. Section \ref{sec:matheur} describes a matheuristic algorithm strongly related to the previous formulations.
Extensive computational results are presented in Section \ref{sec:results}. Section \ref{sec:conclusions} finally concludes the paper.
}

\section{Related literature} \label{sec:lit}
{The Flying Sidekick Traveling Salesman Problem (FSTSP) is a generalization of the (TSP) and the {\em vehicle routing problem} (VRP), and thus is a NP-hard problem. The amount of literature that considers optimization problems related to drones or trucks and drones is limited but it has been experiencing a boom in the last few years.}

 We start our review by considering the problems that Otto et al. \cite{ottooptimization} classify as {\em Drones and vehicles performing independent tasks}. Murray and Chu \cite{murray2015flying} propose the {\em parallel drone scheduling TSP} (PDSTSP), where a fleet of drones can serve customers only departing from the depot. The remaining customers are served by a truck. They propose MILP formulations and simple greedy heuristics for both problems. Mbiadou Saleu et al. \cite{mbiadou2018iterative} propose a two step iterative heuristic based on dynamic programming for the same problem, while Dell'Amico at al. \cite{noi} propose some matheuristic approaches.
Ulmer and Thomas \cite{ulmer2018same} study a dynamic variant of the PDSTSP called the {\em Same-day delivery with heterogeneous fleets of drones and vehicles}, where requests arrive dynamically and they need to be allocated to drones or truck maximizing the number of  served customers. They solve the problem with an approximate dynamic programming known as parametric policy function approximation.
Sawadsitang et al. \cite{sawadsitang2018supplier} propose the {\em joint ground and aerial package delivery service}, that can be considered as a PDSTSP with uncertainty and multiple trucks that they represent as a three stage stochastic model. They solve the model with an L-shaped method.
%
%\subsection{{Trucks and drones working together}}
The following papers consider the coupled interventions of drones and trucks. Savuran and Karakaya \cite{savuran2015route} study a problem where a truck follows a linear path, in the meanwhile a drone is launched from the truck and must return to the truck after performing a, let's say, open TSP to visit all the targets. They solve the problem with a genetic algorithm.  Boysen et al. \cite{boysen2018drone} consider a fixed truck route where the truck represents a loading platform for the drones. The truck and the drone can wait for each other.
Mourelo Ferrandez et al. \cite{ferrandez2016optimization} and Chang and Lee \cite{chang2018optimal}  independently propose works that couple drones and a truck. A set of delivery customers is given, and these custmers are clusterized by using a K-means method. Thus a TSP is solved among the centroid of each cluster that are the points where the truck stops to launch one or more drones. The TSP part is solved with a genetic algorithm.
{Campbell et al. \cite{campbell2017strategic} and Carlsson and Song \cite{carlsson2017coordinated}  treat problems with trucks and drones where, differently from other papers, the demand is continuously distributed.}
{A two-echelon related problem is the one solved by Boysen et al. \cite{boysen2018scheduling} that considers terrestrial robots coupled with a truck. The truck can follow a path where only the starting point is decided and can travel among depots, to collect new robots, or drop-off points, to leave robots. Robots can leave the truck, serve customers, and return to one of the depots. The objective is to minimize the weighted number of late deliveries. They propose MILP model and a multi start local search algorithm. }

%\subsection{{Trucks and drones synchronized}}
We consider, now, routing problems where trucks are equipped with drones, both vehicles can be used to deliver packages to customers, {and synchronization is important}. Otto et al. \cite{ottooptimization} classify these problems under the name {\em Drones and vehicles as synchronized working units}.

\subsection{The Flying sidekick TSP and the TSP with drone}

Murray and Chu \cite{murray2015flying} define and study the {\em Flying sidekick traveling salesman problem} (FSTSP). In the FSTSP the truck and the drone can cooperate to serve customers.
In this case the drone starts from its vehicle in a vertex of the network (launch), performs a delivery to a customer, and returns to the truck (rendezvous) in a vertex of the network. The truck and the drone must be synchronized and thus wait for each other. The objective function is normally to minimize the completion time.
 Customers can be visited only once, and some customers can be visited only by the truck because their request cannot be fulfilled by the drone. The drone cannot return at the launching point.
Agatz et al. \cite{agatz2018optimization} solve the {\em TSP with Drone} (TSP-D). In this case each customer has to be visited at least once by one of the vehicles, but they can be visited more than once by the truck if it is convenient for drone launching and return. Launching and rendezvous can coincide. Some vertices cannot be visited by drones. Endurance is unlimited and launching and rendezvous times are considered negligible. The authors present an ILP model and propose route first-cluster second heuristics based on local search and dynamic programming. Bouman et al. \cite{bouman2018dynamic} solve the TSP-D with dynamic programming.
Ha et al. \cite{ha2015heuristic} also use the name TSP-D, albeit in this case nodes cannot be visited multiple times and launch and rendezvous of a sortie are not allowed to happen at the same vertex. The authors declare that the considered problem shared the FSTSP characteristics, but they called it TSP-D. They propose two heuristic algorithms: a route first-cluster second one and a cluster first-route second one.
{In \cite{ha2017min}, the same authors solve a similar problem with a different objective function, made of four components, each one with a weight: the total distance traveled by the truck, the one travelled by the drone, and the waiting time of the truck and the drone. They present a MILP based on Murray and Chu's one and two heuristics.
 Ha et al. \cite{ha2018hybrid} propose a genetic algorithm to solve both the minimum time and the minimum cost problems.}
{Liu et al. \cite{liu2018optimization} solved the TSP-D with a genetic algorithm that keeps truck of both the minimum time and the minimum energy solutions.}
Ponza's thesis \cite{ponza2016optimization} tackles the FSTSP proposing a modified MILP formulation with respect to the Murray and Chu's \cite{murray2015flying} one and solve it with a simulated annealing algorithm. The new formulation, with respect to the FSTPS one, among the other few differences, does not allow the drone to wait at customers nodes.
{de Freitas Penna \cite{de2018randomized} propose a randomized variable neighborhood descent (VND) for the FSTSP that starts form a TSP solved with Concorde. Some customers are then allocated to the drones. Afterwards they apply five neighborhood in a VND framework. In
\cite{freitas2018variable}, the same authors propose a general Variable Neighborhood Search algorithm for the FSTSP that is a similar to the previous one. }
{Poikonen et al. \cite{poikonen2019branch} propose a branch-and-bound for the TSP-D.}
{Yurek and Ozmutlu \cite{yurek2018decomposition} propose a decomposition-based iterative heuristic for the TSP-D. They work on instances with up to 20 customers, with the exact methods able to solve instances up to 12 customers.}
{Jeong et al. \cite{jeong2019truck} solve the {FSTSP with energy consumption and no fly zones.} This problem accounts for the parcel weight on drone energy consumption and the drone cannot fly over some restricted flying areas. They propose a MILP formulation based on the Murray and Chu one and an evolutionary based heuristic.}
{Marinelli et al. \cite{marinelli2017route} study the en-route TSP-D, where drone can start and return on the arcs travelled by the truck. The launch and return points are decided so to minimize the waiting times. They solve the proposed problem with a GRASP based on the method by Ha et al. \cite{ha2017min} where the initial solution is obtained with the Lin-Kernighan algorithm.}

\subsection{Multiple trucks and drones}
Wang et al. \cite{wang2017vehicle} define the {\em vehicle routing problem with drones} (VRPD), where a homogeneous fleet of trucks equipped with a not necessary unitary number of drones delivers parcels to customers. Drones can be launched from trucks at depot or at any customer vertex. Each drone must return to the same truck and also at the same node where it has been launched.
Poikonen et al. \cite{poikonen2017vehicle} extend the worst-case results considering different metrics for trucks and drones, considering limited drone batteries, and evaluating different objective functions.
Daknama and Kraus \cite{daknama2017vehicle} present and solve the {\em Vehicle routing with drones} where multiple vehicles and drones can be used for deliveries. They minimize the average delivery times instead of the completion time. No mathematical model is presented, but they solve the problem by first solving a multiple TSP heuristically and then introducing drones. Local search procedures are thus applied to improve the solution.

\section{Problem Description} \label{sec:prob}
The FSTSP, firstly defined by Murray and Chu \cite{murray2015flying}, is the problem of serving a set of customers $C = \{1,\dots,c\}$ with either a truck or a drone. The truck starts from the depot 0 and returns to the final depot $c+1$, and is equipped with a flying drone that can be used to serve one customer at a time, in parallel to the truck. A drone service is called {\em sortie}, defined by a launching node, a served customer, and a rendezvous node. All customers of $C$ can be served by the truck, but only a subset $C' \subseteq C$ can be served by the drone with a sortie. The problem is built on digraph $G=(N,A)$, where the set $N = \{0,1,\dots,c+1\}$ represents all the nodes, while we define $N_0 = \{0,1,\dots,c\}$ and $N_+=\{1,\dots,c+1\}$. Let $A$ be the set of all the arcs $(i,j), i \in N_0, j \in N_+$, $i\neq j$. Each arc $(i,j)$ is associated with two non-negatives traveling times: $\tau_{ij}^T$ and $\tau_{ij}^D$, that represent the time for traveling that arc by the truck and by the drone, respectively. The travel time matrices of the drone and the truck are normally different. 
Nodes 0 and $c+1$ represent the same physical point, the depot, and the traveling time between them is set to 0.
Serving times at customers for both drone and truck are included in the travel times, while the time for preparing the drone at launch is given by $\sigma^L$ and  the rendezvous time is given by $\sigma^R$. No launch time is considered when the sortie starts form the depot.
The drone have a battery limit (endurance) of $E$ time units, that constraints its use. {Rendezvous time $\sigma^R$ contributes to the endurance computation while $\sigma^L$ does not, since the drone lies on the truck when it is prepared for the launch.}

A {sortie} is formally defined by a triplet $\langle i,j,k \rangle$, ($i\ne j\ne k$) where $i \in N_0$ is the launching node, $j \in C'$ the customer to serve, and $k \in N_+ $ the rendezvous node.
Let $F$ be the set of all sorties that can be performed within the endurance time $E$ ($\tau_{ij}^D + \tau_{jk}^D + \sigma^R \le E$).

The drone can be launched from the truck only when the truck is stopped at a customer or at the depot;
the drone cannot leave the depot before the truck starts its route.
The truck can keep serving customers while the drone is performing a sortie.  A synchronization is required: the vehicle (drone or truck) that arrives first at a rendezvous point has to wait for the other.
The objective of the optimization is to minimize the completion time, that is the moment when the last vehicle arrives at the depot.

\section{Formulations}\label{sec:form}
In the following we propose three new formulations for the FSTSP. In the first model sorties are represented as 3-indexed variables, in the last two as two-indexed variables. {The novelty of these models is that we use only one set of time variables with respect to formerly presented formulations, that used one set of time variables for the truck and one for the drone. We thus halved the number of time variables. The new time variables are used to model synchronization and compute waiting times. They are not needed in the objective function and thus we can write more compact models. Another novelty is that in the first two formulations we include a binary variable that states the position of the drone with respect to the truck. This helps in avoiding crossing sorties (see below). In the third formulation we manage crossing sorties by means of a set of constraints included in the model.}

\subsection{{A 3-indexed Formulation}}%
The first formulation we present,  {3IF}, is built on the formulation DMN proposed by Dell'Amico et al. \cite{dell2019drone}, where the truck route is represented
making use of the variable $x_{ij} =1$ if the node $j \in N_+$ is visited after node $i \in N_0, j \ne i$, and 0 otherwise.
The drone sorties are represented by a three indexed variable $y_{ijk}, \langle i,j,k \rangle \in F$, that equals 1 if the sortie is performed, 0 otherwise. Non-negative variables $w_i$ represent the time that the truck waits for the drone at node $i \in N$. {Non-negative variables $t_i, i \in N$ are used to represent the time synchronization and identifies one of the novelties of this model. Another novelty is the variable $z_i$ that equals 1 if the drone is on the truck at node $i \in N$, 0 otherwise.}
In the following we describe the formulation 3IF step by step.\\

\noindent \emph{Objective function}\\
The objective function \eqref{eq:fobb_dmn1} aims at minimizing the arrival of truck and drone at the final depot. {We clarify that the drone can arrive at the depot after the truck and this is described by the waiting time variable $w_{c+1}$ }. The completion time can be decomposed into the truck route traveling time, the time needed for launching and collecting a drone, and the time the truck waits for the drone.

\begin{align}
    &\min \sum_{(i,j) \in A}\tau_{ij}^Tx_{ij} + \sigma^R\sum_{\langle 0,j,k\rangle \in F}y_{0jk} +  (\sigma^L + \sigma^R)\sum_{\langle i,j,k\rangle \in F, i \neq 0}y_{ijk}  + \sum_{i \in N_+}w_i \label{eq:fobb_dmn1}
\end{align}

\noindent \emph{Customer covering}\\
Constraints \eqref{eq:in_dmn1} and \eqref{eq:out_dmn1}  enforce one of the two vehicles to serve each customer exactly once.
\begin{align}
      &\sum_{i|(i,j)\in A}x_{ij} + \sum_{i,k |\langle i,j,k\rangle \in F}y_{ijk}  = 1 \quad j \in C  \label{eq:in_dmn1}\\
     &\sum_{i|(j,i)\in A}x_{ji} + \sum_{i,k |\langle i,j,k\rangle \in F}y_{ijk}  = 1 \quad j \in C  \label{eq:out_dmn1}
\end{align}

\noindent \emph{Truck routing constraints}\\
In \eqref{eq:startend_dmn1} we impose that the truck starts and finishes its journey at the depots.
\begin{align}
    &\sum_{j \in N_+}x_{0j} = \sum_{i \in N_0}x_{i,c+1} =1 \label{eq:startend_dmn1}
\end{align}
Note that the combination of \eqref{eq:in_dmn1} and \eqref{eq:out_dmn1} also implies the flow conservation for the truck, i.e., 
\begin{align}
&\sum_{i | (i,j) \in A}x_{ij}=\sum_{i| (j,i)\in A}x_{ji} \quad j \in C \label{eq:flow_dmn1}
\end{align}

\noindent \emph{Timing constraints}\\
 In order to guarantee the timing constraints, we impose constraints \eqref{eq:time1_dmn1} to ensure that, if arc $(i,j)$ is traveled by the truck, then the time in $j$ is at least the time in $i$ plus the time needed for traveling the arc. In \eqref{eq:time2_dnm1} we assure that if there is a sortie $\langle i,k,j \rangle$ then the time $t_j$ should be at least $t_i$ plus the time for performing the sortie.
If the drone arrives after the truck at the rendezvous point, then the truck must wait at least the difference between the arrival time of the drone and the arrival time of the truck in that node, that is imposed by constraint \eqref{eq:timewait_dmn1}. Note that if the drone arrives before the truck, it waits while flying, this waiting time is absorbed by time $t_j$ of constraint \eqref{eq:time1_dmn1}, which is included in the truck route in the objective function.
\begin{align}
     &t_j \ge t_i +\tau_{ij}^T - M(1 - x_{ij}) \quad (i,j) \in A \label{eq:time1_dmn1}\\
    &t_j \ge t_i + \tau_{ik}^D + \tau_{kj}^D - M( 1- \sum_{j | \langle i,k,j\rangle \in F}y_{ikj}) \quad (i,j) \in A \label{eq:time2_dnm1}\\
    &w_j \ge t_j - t_i - \tau_{ij}^T  - M(1-x_{ij}) \quad (i,j) \in A \label{eq:timewait_dmn1}
\end{align}
Note that the launching and rendezvous service times are not considered in those constraints because it is not necessary to include them in our variables $t$. Indeed, launching and rendezvous service times should be included in both the truck and drone timing constraints \eqref{eq:time1_dmn1} and \eqref{eq:time2_dnm1}, however they have no effect in the computation of the waiting times \eqref{eq:timewait_dmn1} since they appear as a constant in both the truck and drone time.
Constraints \eqref{eq:time1_dmn1} and \eqref{eq:time2_dnm1} also have the side effect of avoiding backward sorties, i.e. sorties $\langle i,k,j \rangle \in F$ such that node $j$ in visited by the truck before node $i$.
In contrast with the other models in the literature, $t$ variables do not represent the exact time of visiting one node, but they allow the model to respect the truck-drone timing, to compute the waiting times, and allow us to obtain a more compact model.\\

\noindent \emph{Drone battery endurance constraints}\\
In next constraint \eqref{eq:endurance_dmn1} we assure that if a sortie $\langle i,j,k \rangle$ is performed, the elapsed time from the launch to the rendezvous respects the drone endurance $E$. {The rendezvous service time $\sigma^R$ is included in the drone time.}%
\begin{align}
    &t_k - t_i + \sigma^R - M(1-\sum_{j | \langle i,j,k\rangle \in F}y_{ijk}) \le E \quad (i,k) \in A \label{eq:endurance_dmn1}
\end{align}

\noindent \emph{x-z and y-z linking constraints}\\
Constraints \eqref{eq:zeta2_dmn1} state that a drone can be on the truck at node $i$ if the truck enters node $i$ {or, in turn, using \eqref{eq:flow_dmn1}, that the truck exits node $i$.} Constraints \eqref{eq:zeta3_dmn1} state that a sortie can start in $i$ if  variable $z_i$ equals 1. Those constraints, in addition with constraints \eqref{eq:zeta1_dmn1} are also used to avoid crossing sorties.
A pair of sorties $\langle i,j,k \rangle, \langle i',j',k'\rangle \in F$ are called \textit{crossing} if $i$ is visited by the truck before $i'$, but $k$ is not visited before $i'$. In other words, the second sortie starts before the first one is terminated.
Constraints \eqref{eq:zeta1_dmn1} regulate $z$ variables, modeling the presence of the drone over the truck route. 
{Consider the first node $i$ on the truck route from which starts a sortie $\langle i, k, l \rangle \in F$ and arc $(i,j) \in A$ such that $x_{ij} = 1$. No sortie can arrive in node $j$ (apart from the one starting in node $i$) since \eqref{eq:time1_dmn1} and \eqref{eq:time2_dnm1} forbid backward sorties.
Hence, given $z_i =1$ because of \eqref{eq:zeta3_dmn1}, \eqref{eq:zeta1_dmn1} impose $z_j = 0$, while if the sortie that starts in $i$ returns in $j$, $z_j$ can take value 1, which means that a new sortie can start in $j$. Coming back to the general case, where the sortie do not return in $j$, the $z$ variables along the path after the sortie launch in $i$ and before the rendezvous in $l$ remain to value 0, imposing that no other sortie can start until the truck reaches node $l$ where the drone returns. After that, the $z$ variables can take value 1, which means that a new sortie can start.}
\begin{align}
    &z_i \le \sum_{j| (j,i) \in A}x_{ji}\quad i \in N_+ \label{eq:zeta2_dmn1}\\
    &\sum_{j, k | \langle i,j,k\rangle \in F}y_{ijk} \le z_i \quad i \in N_0\label{eq:zeta3_dmn1}\\
        &z_j \le z_i - x_{ij} + \sum_{l,k|\langle l,k,j \rangle \in F}y_{lkj}- \sum_{k,l|\langle i,k,l \rangle \in F}y_{ikl} +1 \quad (i,j) \in A\label{eq:zeta1_dmn1}
\end{align}

\noindent \emph{Variable bounds}
\begin{align}
    &t_0 = 0 \label{eq:tzero_dmn1} \\
    &t_i, w_i \in \mathbb{R}^+ \quad i \in N_+ \label{eq:pos_dmn1}\\
    &z_i \in \{0,1\}\quad  i \in N \label{eq:binz_dmn1}\\
    &x_{ij} \in \{0,1\} \quad (i,j) \in A \label{eq:binz_dmn1}\\
    &y_{ijk} \in \{0,1\} \quad \langle i,j,k\rangle \in F
\end{align}

\subsection{{A 2-indexed Formulation}} \label{2if}
\begin{sloppypar}
Formulation {2IF} is based on the formulation DMN2 presented in Dell'Amico et al. \cite{dell2019drone}. Similarly to the previous model, the main difference with DMN2 is that 2IF makes use of only one set of time variables to represent synchronization, halving the number of time variables and diminishing the number of  `big-M' constraints required. The other main difference is the use of a binary variable $z$ representing the presence or absence of the drone on the truck along the route.
We model the truck route as in the 3IF; however, in this formulation we model the sorties with a pair of two-indexed binary variables: one for the launch and one for the rendezvous. Variable $\overrightarrow{g}_{ij}$ takes value 1 if the drone is launched in $i \in N_0$ and serves the customer $j \in C'$ and $\overleftarrow{g}_{jk}$ is 1 if the drone returns to node $k \in N_+$ after visiting customer $j \in C'$.
To reduce the variables we preliminary fix to zero all  those corresponding  to arcs with flying time exceeding the battery limit, i.e., we set $\overrightarrow{g}_{ij}=0$ for all $ (i,j)\in A: \tau^D_{ij}>E$ and
$\overleftarrow{g}_{jk}=0$  for all  $(j,k)\in A: \tau^D_{jk}+\sigma_R>E$. We also fix to zero variables that do not allow to complete a feasible drone fly: $\overrightarrow{g}_{ij}=0$, $(i,j)\in A, j\not\in C'$; $\overleftarrow{g}_{jk}=0$, $(j,k)\in A, j\not\in C'$; $\overrightarrow{g}_{ic+1}=0$ $i\in N$
 and $\overleftarrow{g}_{j0}=0$ $j\in N$.
Variables $x, t, w, z$ are the same as for the previous model.\\
\end{sloppypar}
{In formulation 2IF truck routing constraints are the same as for \eqref{eq:startend_dmn1} and \eqref{eq:flow_dmn1}. We also import  form 3IF the timing constraints \eqref{eq:time1_dmn1} and \eqref{eq:timewait_dmn1}, constraint \eqref{eq:zeta2_dmn1}, and variable bounds \eqref{eq:tzero_dmn1}--\eqref{eq:binz_dmn1}. The remaining }components are as follows.\\

\noindent \emph{Objective function}\\
The objective function \eqref{eq:dmn4_fobb} {is the adjustment of \eqref{eq:fobb_dmn1} to the new sets of variables.}
\begin{align}
 &\min  \sum_{(i,j)\in A}\tau_{ij}^Tx_{ij} + \sigma^L\sum_{\substack{(i,j)\in A\\i\neq 0}}\overrightarrow{g}_{ij} + \sigma^R\sum_{(j,k)\in A}\overleftarrow{g}_{jk} + \sum_{i \in N_+}w_i\label{eq:dmn4_fobb}
\end{align}

\noindent \emph{Customer covering} \\
Constraints \eqref{eq:dmn4_ass} and  \eqref{eq:dmn4_ass2} impose that all customers must be served either by the truck or by the drone:
\begin{align}
&\sum_{i|(i,j)\in A} x_{ij}  +  \sum_{i|(i,j)\in A} \overrightarrow{g}_{ij} = 1 \quad  j \in C\label{eq:dmn4_ass}\\
&\sum_{i|(j,i)\in A} x_{ji}  +  \sum_{i|(j,i)\in A} \overleftarrow{g}_{ji} = 1 \quad  j \in C\label{eq:dmn4_ass2}
\end{align}

\noindent \emph{Timing constraints}\\
In addition to constraints \eqref{eq:time1_dmn1} and \eqref{eq:timewait_dmn1}, we impose constraints \eqref{eq:time2_dmn4} and \eqref{eq:time3_dmn4} to update the times of the two components of a sortie (launch and rendezvous) separately.
\begin{align}
    &t_j \ge t_i +\tau_{ij}^D  - M(1-\overrightarrow{g}_{ij}) \quad (i, j) \in A \label{eq:time2_dmn4}\\
    &t_k \ge t_j + \tau_{jk}^D - M(1 - \overleftarrow{g}_{jk}) \quad (j,k) \in A \label{eq:time3_dmn4}
\end{align}

\noindent \emph{Drone battery endurance constraint}\\
If a sortie is performed then its total time should respect the battery endurance, as for constraints \eqref{eq:dmn5_energy}.
\begin{align}
    &t_k - t_i  +  \sigma^R - M(2-\overrightarrow{g}_{ij}- \overleftarrow{g}_{jk})\le E \quad i \in N_0, j \in C', k \in N_+ \label{eq:dmn5_energy}
\end{align}

\noindent \emph{Sortie congruence}\\
Constraints \eqref{eq:dmn5_flowg} impose that if a drone serves a customer then it must also return to the truck.
\begin{align}
      & \sum_{i|(i,j)\in A} \overrightarrow{g}_{ij} = \sum_{k|(j,k)\in A} \overleftarrow{g}_{jk} \quad j \in C' \label{eq:dmn5_flowg}
\end{align}

\noindent \emph{x-z and g-z linking constraints}\\
Similarly to previous formulation and in addition to constraints \eqref{eq:zeta2_dmn1}, we use \eqref{eq:zeta3_dmn5} to impose that a sortie can start in node $i$ only if $z_i= 1$ and \eqref{eq:zeta1_dmn5} to avoid crossing sorties:
\begin{align}
    &\sum_{j \in C'}\overrightarrow{g}_{ij} \le z_i \quad i \in N_0\label{eq:zeta3_dmn5}\\
    &z_j \le z_i - x_{ij} + \sum_{k \in C'}(\overleftarrow{g}_{kj}- \overrightarrow{g}_{ik}) +1 \quad (i,j) \in A\label{eq:zeta1_dmn5}
\end{align}

\noindent \emph{Variable bounds}\\
In addition to
\eqref{eq:tzero_dmn1}-\eqref{eq:binz_dmn1}, we introduce the constraint on variables $g$:
\begin{align}
    &\overrightarrow{g}_{ij}, \overleftarrow{g}_{ij} \in \{0,1\} \quad (i,j) \in A \label{eq:binz_dmn5}
\end{align}

\subsubsection{Inequalities}

The following inequalities \eqref{eq:ineq_twoindex1} are a relaxed version of constraint \eqref{eq:dmn5_energy}, they avoid sorties longer than the maximum endurance. Notwithstanding this version of the constraint does not take into account the drone waiting time in case it arrives to the rendezvous point before the truck, it has a strong impact on the solution time of the model, since it is not affected by large constants ($M$).
\begin{align}
    \sum_{i \in N_0}\tau_{ij}^D\overrightarrow{g}_{ij} + \sum_{k \in N_+}\tau^D_{jk}\overleftarrow{g}_{jk} + \sigma^R\le E \quad j \in C'
    \label{eq:ineq_twoindex1}
\end{align}
These constraints are used for the experiments reported in Section \ref{sec:results}.

\subsection{A modified 2-indexed Formulation} \label{sec:dmn3}

{In this section we present a mathematical formulation built upon the formulation {2IF} with the difference that variables $z$ are not necessary anymore to avoid infeasible sorties. For doing so we make use of the crossing sorties elimination constraints proposed by Dell'Amico et al. \cite{dell2019drone}. These constraints are exponentially many and we separate them in a branch-and-cut (BC) fashion. We refer to this formulation as {2IF-BC}.}

The set of  variables used is the same as in 2IF, with the exclusion of variables $z$. The following constraints are added to \eqref{eq:startend_dmn1}-\eqref{eq:time1_dmn1}, \eqref{eq:timewait_dmn1}, \eqref{eq:tzero_dmn1},\eqref{eq:pos_dmn1},\eqref{eq:binz_dmn1}, \eqref{eq:dmn4_fobb}-\eqref{eq:dmn5_flowg}, \eqref{eq:binz_dmn5}:\\

\noindent\emph{$x$-$g$ coupling constraints}\\
Since this formulation does not make use of $z$ variables,  we need to impose that a sortie can start and end in one node only if the truck is there, respectively:
\begin{align}
& \sum_{j|(i,j)\in A}\overrightarrow{g}_{ij} \le \sum_{h|(i,h)\in A}x_{ih} \quad i \in  N_0 \label{eq:dmn5_link}\\
& \sum_{i|(i,j)\in A}\overleftarrow{g}_{ij} \le \sum_{h|(h,j)\in A}x_{hj} \quad j \in  N_+ \label{eq:dmn5_link2}
\end{align}

\noindent\emph{Crossing Sorties Elimination Constraints}\\
To avoid crossing sorties {we need to impose} the inequalities \eqref{eq:dmn5_CSEC_2}, {that we report in its tournament version hereafter (see \cite{dell2019drone} for details)}.
Let $i \in N_0,l\in C$ be the starting and ending {vertices} of the truck path $P$ from $i$ to $l$ and $P=\{v(1),v(2),\dots,v(q)\}$ with $v(1)=i, v(q)=l$. We assume that exist two sorties defined by $\overrightarrow{g}_{ij}>0$ and $\overrightarrow{g}_{lm}>0$ and that there is no node $k\in P\setminus\{i,l\}$ with $\overleftarrow{g}_{jk}>0$. In this case the second sortie starts before the first sortie is terminated an the following ``tournament" crossing sorties elimination holds:
\allowdisplaybreaks
\begin{align}
& \sum_{h = 1}^{|P| - 1}\sum_{j  = h+1}^{|P|} x_{v(h)v(j)} + \sum_{\substack{(i,j)\in A, \\j\not\in P}}\overrightarrow{g}_{ij} +\sum_{\substack{(l,j)\in A, \\j\not\in P}}\overrightarrow{g}_{lj}
 \le |P|& P \in \mathcal{P}\label{eq:dmn5_CSEC_2}
\end{align}
where $\mathcal{P}$ defines the set of all the paths with the described characteristics. \\

\noindent \emph{{Trivial infeasible sorties avoidance}}\\
We finally impose the following constraints to avoid infeasibilities:
\begin{align}
&\overrightarrow{g}_{ij} + \overleftarrow{g}_{ij} \le 1 \quad (i,j) \in A\\
&\overrightarrow{g}_{ij} + \overleftarrow{g}_{ji} \le 1 \quad (i,j) \in A
\end{align}

\section{{A Random Restart Local Search matheuristic algorithm}}\label{sec:matheur}
\begin{sloppypar}
In this section we propose a matheuristic Random Restart Local Search (RRLS) algorithm (see, e.g., \cite{lourencco2019iterated}) {that} takes advantage of the abilities of the MILP solver to provide high quality heuristic solutions, once a truck tour has been identified, and under some super-imposed precedence constraints among customers.
\end{sloppypar}

The idea is to optimize the truck tour with state-of-the-art heuristics as a classic TSP, and to delegate the MILP model 2IF to insert drone missions into such a truck tour, with the additional constraint that the visiting of the customers has to fulfil the order imposed by the TSP solution previously calculated. 

Once a FSTSP solution is provided, an iterative mechanism is entered, where the new truck tour returned by the MILP solver (and typically not covering all the customers) is reoptimized by a TSP solver and passed again to the MILP. Once a local minimum is reached, the algorithm restarts from a new TSP solution calculated by adding some random noise to the distances, in order to increase the exploration of the FSTSP search space. The algorithm stops when a given computation time has elapsed.

The method can be formally described through the following pseudocode:
\begin{enumerate}
\item $Bestcost=+\infty$
\item A TSP instance with the customers of the FSTSP is solved with the algorithm LKH \cite{hel06}, obtaining a sequence of customers $s$ such that customer $s_i$ is in position $i$ in the sequence. (Note that at this stage the tour $s$ covers all the customers, while in the next steps this might not be always the case.)
\item The MILP model 2IF described in Section \ref{2if} is then solved with the following additional constraints:
 \begin{equation}
 t_{s_i} \leq t_{s_{i+1}} \ \ \forall i \in N_0
 \end{equation}
The computation is eventually truncated after 30 seconds (or once a feasible solution is retrieved if this time does not suffice), and a solution $Sol$ with cost $c(Sol)$ is stored.\\
 Let $s'$ denote the sequence of customers visited by the truck in $Sol$.
\item If $c(Sol) < Bestcost$ then $Bestcost = c(Sol)$ and $Bestsol = Sol$.
\item Algorithm LKH \cite{hel06} is run on the customers contained in the sequence $s'$ to improve the truck tour, obtaining the optimized sequence $s$. 
\item
If $s \equiv s'$ then
 a new TSP route is generated by running algorithm LKH \cite{hel06} on an temporary graph with artificial travel times $\overline{t_{ij}^T} = (1+ rand(0,0.5)) t_{ij}^T$, where $rand(0, 0.5)$ is a random number between 0 and 0.5. The idea is to restart the local search from a different area of the search space.
\item If the exit criterion is not met (this is typically a maximum computation time), go to step 3.
\end{enumerate}

\section{Computational Experiments}\label{sec:results}
The algorithms have been implemented in ANSI C. The mathematical models have been solved by Gurobi 8.1 on an Intel Core i3-2100 CPU, with 3.10 GHz and 8.00 GB
of RAM.
Formulation {3IF and 2IF} have been solved directly by the solver, while {2IF-BC} required a branch-and-cut implementation, since it includes {constraints \eqref{eq:dmn5_CSEC_2}, that are exponentially
many}. To separate these constraints, we have considered the residual graph $G' = (N, A')$ obtained from $G$ by selecting the only arcs associated with a non zero variable ($x$, $\overrightarrow{g}$, $\overleftarrow{g}$) in the continuous relaxation of the model. For the crossing sorties elimination constraints we explore the graph starting from depot 0, until a truck path violating one of the constraints is identified, if any. {After preliminary computational tests, we have observed that formulation 2IF-BC benefited from separating these constraints only for integer solutions. In such a case the overall procedure has a time complexity $O(|A'|)$.}
The matheuristic algorithm has been solved on a Intel(R) Xeon(R) E5-2620 v4 2.10 GHz. All computation times are obtained with single thread runs.  {The CPU marks on single thread based on the results in \url{https://www.cpubenchmark.net/} for the two used machines are very similar, being 1.584 and 1.547, respectively.}
 
\subsection{Small Instances}\label{sec:results_small}
We first run our methods on the 72 randomly generated benchmark instances proposed by Murray and Chu \cite{murray2015flying}, with 10 customers and the endurance $E$ set to either 20 or 40.
In Table \ref{tab:smallinst20} one can find the results of our methods for $E=20$, and in Table \ref{tab:smallinst40} the results for $E=40$.
In the Tables we report the instance name, the optimal solution value {(for instances `37v5', `37v9', `37v10', and `43v9' with $E=20$ and for instances `37v1', `37v5', `37v6', `37v8', `37v9', `37v10', `40v4', `43v1',  `43v5',  `43v6', and  `43v9' with $E=40$ documented here for the first time)}, and the time needed by each of the exact methods to certify a solution. We also report the gap, calculated as $gap\% = 100 \cdot (cost(RRLS) - opt)/opt$, of the solutions provided by the heuristic RRLS, for which a maximum computation time of 20 seconds was considered.

\begin{table}[htbp]
  \centering
  \footnotesize
  \caption{Results on small  instances from \cite{murray2015flying} with endurance $E=20$.}
    \begin{tabular}{crrrrc}
    \toprule
    Instance  & \multicolumn{1}{l}{opt} & \multicolumn{1}{c}{3IF} & \multicolumn{1}{c}{2IF} & \multicolumn{1}{c}{2IF-BC} & \multicolumn{1}{c}{RRLS (20s)} \\
     & & \multicolumn{1}{c}{sec} & \multicolumn{1}{c}{sec} & \multicolumn{1}{c}{sec} & \multicolumn{1}{c}{gap\%} \\
    \cmidrule(lr){1-1}  \cmidrule(lr){2-2} \cmidrule(lr){6-6} \cmidrule(lr){3-5}
37v1	&	57.45	&	16.45	&	11.86	&	1.23	&	0.00	\\
37v2	&	53.79	&	5.83	&	4.18	&	0.79	&	0.00	\\
37v3	&	54.66	&	6.81	&	7.60	&	6.21	&	0.00	\\
37v4	&	67.46	&	5.09	&	2.28	&	5.09	&	0.00	\\
37v5	&	51.78	&	6114.57	&	2254.65	&	406.43	&	1.28	\\
37v6	&	48.60	&	1786.79	&	672.74	&	128.31	&	0.00	\\
37v7	&	49.58	&	748.76	&	113.76	&	18.38	&	0.00	\\
37v8	&	62.38	&	324.32	&	142.60	&	43.12	&	0.00	\\
37v9	&	43.48	&	3215.10	&	1551.05	&	361.49	&	6.68	\\
37v10	&	41.91	&	229.14	&	271.68	&	211.62	&	0.00	\\
37v11	&	42.90	&	17.82	&	31.90	&	37.47	&	2.33	\\
37v12	&	56.85	&	36.04	&	52.49	&	38.31	&	0.00	\\
40v1	&	49.43	&	12.04	&	11.86	&	4.03	&	1.95	\\
40v2	&	51.71	&	20.05	&	18.13	&	12.09	&	4.54	\\
40v3	&	57.10	&	8.11	&	13.15	&	6.90	&	0.00	\\
40v4	&	69.90	&	2.72	&	7.39	&	1.98	&	0.00	\\
40v5	&	45.46	&	475.19	&	173.67	&	132.01	&	2.20	\\
40v6	&	44.51	&	56.19	&	16.52	&	9.48	&	0.00	\\
40v7	&	49.90	&	10.74	&	3.92	&	5.27	&	0.00	\\
40v8	&	62.70	&	10.23	&	4.32	&	5.84	&	0.00	\\
40v9	&	42.53	&	18.33	&	7.60	&	14.47	&	0.00	\\
40v10	&	43.08	&	2.16	&	1.54	&	3.64	&	0.00	\\
40v11	&	49.20	&	1.40	&	1.85	&	1.19	&	0.00	\\
40v12	&	62.00	&	1.97	&	1.78	&	1.12	&	0.00	\\
43v1		&	69.59	&	0.62	&	0.62	&	0.53	&	0.00	\\
43v2		&	72.15	&	0.70	&	0.68	&	0.58	&	0.00	\\
43v3		&	77.34	&	0.64	&	1.00	&	0.43	&	0.00	\\
43v4		&	90.14	&	0.81	&	0.56	&	0.53	&	0.00	\\
43v5		&	58.71	&	2131.85	&	1509.45	&	242.87	&	0.78	\\
43v6		&	59.09	&	874.77	&	394.40	&	69.65	&	0.77	\\
43v7		&	65.52	&	284.10	&	122.68	&	8.98	&	0.00	\\
43v8	&	84.81	&	634.77	&	276.14	&	103.32	&	0.00	\\
43v9		&	46.93	&	1340.57	&	1205.13	&	707.32	&	0.00	\\
43v10		&	47.93	&	121.91	&	158.75	&	72.36	&	0.00	\\
43v11		&	57.38	&	1.73	&	6.63	&	11.91	&	0.00	\\
43v12		&	69.20	&	2.05	&	1.58	&	5.28	&	0.00	\\
    \cmidrule(lr){1-1}  \cmidrule(lr){2-2} \cmidrule(lr){6-6} \cmidrule(lr){3-5}
    Avg. time & & 514.45 & 251.56 & 74.45 & 3.15 \\
    Avg. \%gap & &  &  &  & 0.57 \\
    \bottomrule
    \end{tabular}%
  \label{tab:smallinst20}%
\end{table}%

\begin{table}[htbp]
  \centering
  \footnotesize
  \caption{Results on small  instances from \cite{murray2015flying} with endurance $E=40$.}
    \begin{tabular}{crrrrc}
    \toprule
    Instance  & \multicolumn{1}{l}{opt} & \multicolumn{1}{c}{3IF} & \multicolumn{1}{c}{2IF} & \multicolumn{1}{c}{2IF-BC} & \multicolumn{1}{c}{RRLS (20s)} \\
     & & \multicolumn{1}{c}{sec} & \multicolumn{1}{c}{sec} & \multicolumn{1}{c}{sec} & \multicolumn{1}{c}{gap\%} \\
    \cmidrule(lr){1-1}  \cmidrule(lr){2-2} \cmidrule(lr){6-6} \cmidrule(lr){3-5}
37v1	&	50.57	&	7615.59	&	6765.07	&	2731.91	&	0.72	\\
37v2	&	47.31	&	2695.95	&	517.77	&	284.46	&	0.00	\\
37v3	&	53.69	&	1966.82	&	532.39	&	190.49	&	0.00	\\
37v4	&	67.46	&	3720.39	&	1283.28	&	148.43	&	0.00	\\
37v5	&	45.84	&	2364.48	&	2875.92	&	644.23	&	0.00	\\
37v6	&	44.60	&	1682.59	&	458.63	&	253.15	&	0.00	\\
37v7	&	47.62	&	1008.24	&	522.79	&	267.10	&	0.00	\\
37v8	&	60.42	&	1823.89	&	576.20	&	241.61	&	0.00	\\
37v9	&	42.42	&	1501.51	&	544.74	&	1167.62	&	0.00	\\
37v10	&	41.91	&	371.44	&	221.57	&	406.15	&	1.21	\\
37v11	&	42.90	&	30.07	&	23.92	&	51.09	&	0.00	\\
37v12	&	55.70	&	73.00	&	32.36	&	26.62	&	0.00	\\
40v1	&	46.89	&	1160.87	&	464.25	&	252.66	&	0.00	\\
40v2	&	46.42	&	244.42	&	134.70	&	70.15	&	0.00	\\
40v3	&	53.93	&	1418.20	&	280.28	&	150.67	&	2.59	\\
40v4	&	68.40	&	4054.47	&	1439.61	&	252.36	&	0.74	\\
40v5	&	43.53	&	54.99	&	31.91	&	53.75	&	0.00	\\
40v6	&	44.08	&	28.71	&	30.42	&	40.61	&	0.00	\\
40v7	&	49.23	&	1.80	&	1.38	&	1.13	&	0.00	\\
40v8	&	62.03	&	2.40	&	2.95	&	6.76	&	0.00	\\
40v9	&	42.53	&	26.42	&	13.75	&	5.50	&	0.00	\\
40v10	&	43.08	&	1.97	&	5.83	&	5.53	&	0.00	\\
40v11	&	49.20	&	1.07	&	1.35	&	3.38	&	0.00	\\
40v12	&	62.00	&	2.76	&	1.15	&	1.16	&	0.00	\\
43v1	&	57.01	&	18707.53	&	4032.67	&	2312.62	&	0.00	\\
43v2	&	58.05	&	5575.49	&	1978.44	&	1446.23	&	0.00	\\
43v3	&	69.43	&	4901.81	&	1606.24	&	219.80	&	2.12	\\
43v4	&	83.70	&	11045.26	&	1028.85	&	224.88	&	0.00	\\
43v5	&	52.09	&	19768.66	&	13007.44	&	15702.75	&	5.49	\\
43v6	&	52.33	&	5895.86	&	3009.49	&	1921.06	&	4.79	\\
43v7	&	61.88	&	1249.73	&	245.27	&	452.45	&	0.00	\\
43v8	&	73.73	&	407.65	&	177.02	&	123.31	&	0.32	\\
43v9	&	46.93	&	5737.30	&	4689.05	&	5271.24	&	0.00	\\
43v10	&	47.93	&	501.27	&	311.91	&	169.56	&	0.00	\\
43v11	&	56.40	&	1.70	&	1.60	&	11.37	&	0.00	\\
43v12	&	69.20	&	11.58	&	1.04	&	4.21	&	0.00	\\
    \cmidrule(lr){1-1}  \cmidrule(lr){2-2} \cmidrule(lr){6-6} \cmidrule(lr){3-5}
    Avg. time & & 2934.89 & 1301.42 & 975.44 & 2.96 \\
    Avg. \%gap & &  &  &  & 0.50 \\
    \bottomrule
    \end{tabular}%
  \label{tab:smallinst40}%
\end{table}%

The results of Table \ref{tab:smallinst20} indicate a clear ranking of the formulations, with 2IF-BC able to converge faster than the other models, followed by 2IF and 3IF. The RRLS algorithm provides high quality results (gap always below 5\% with a majority of optimal solutions retrieved), with a computation time that is drastically shorter than the exact methods. The exam of Table  \ref{tab:smallinst40} substantially supports the same conclusions, although here formulation 2IF converges faster than 2IF-BC for some of the instances. This can be intuitively explained by considering that with $E=40$ the drone has more freedom, and therefore more cuts need to be iteratively separated.

\begin{figure}
     \centering
     \includegraphics[width=0.6\textwidth]{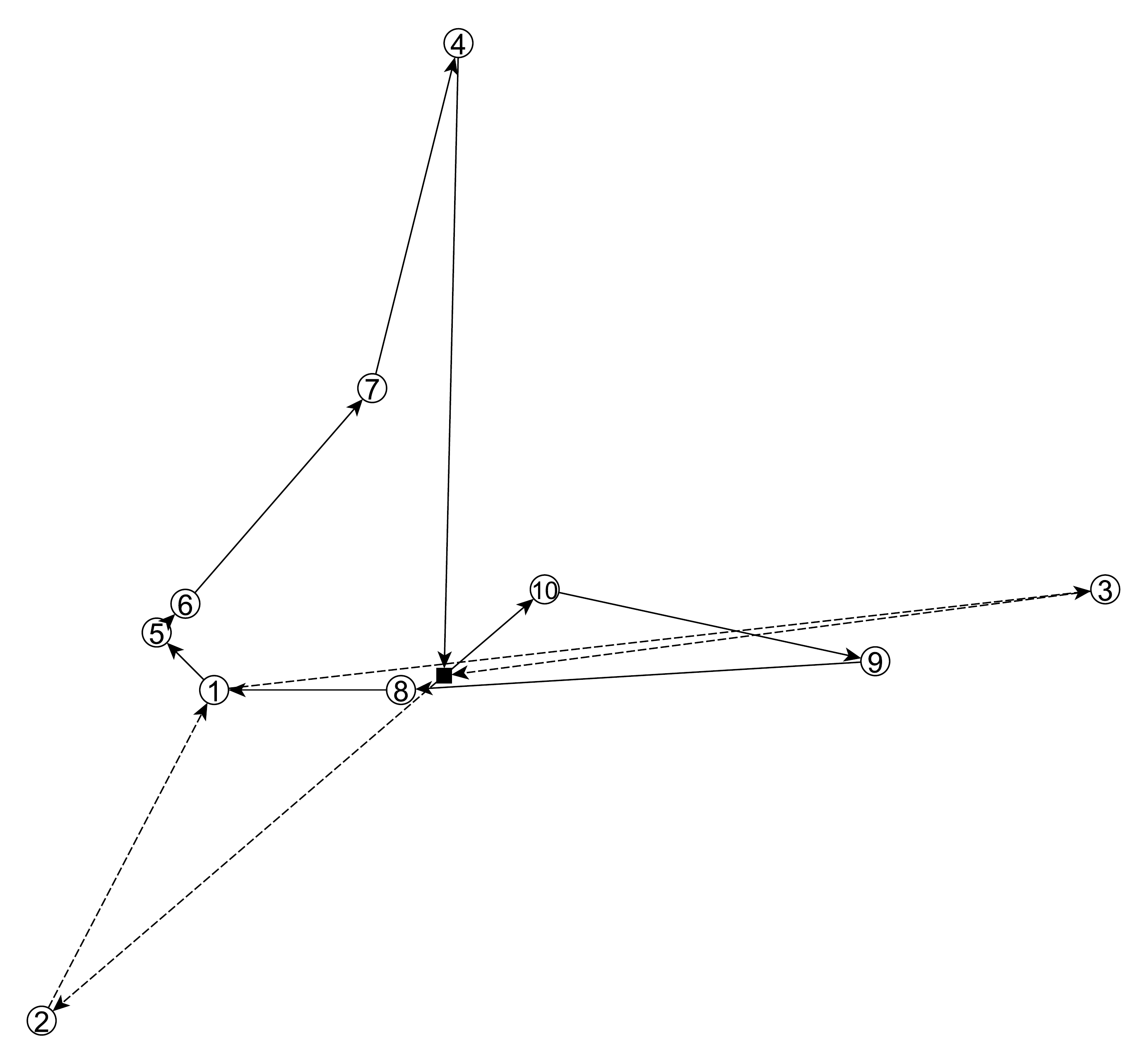}
     \caption{Optimal solution for instance 
     %`20140810T123443v5'
     `43v5' with $E=40$. The black square is the depot, the other nodes are the customers, the solid arcs represent the truck route, while the dashed ones represent the sorties.}
     \label{fig:outlier}
 \end{figure}

\begin{sloppypar}
{It is worth observing that instance  %`20140810T123443v5' 
`43v5' with $E=40$ appears to be an outlier, requiring higher solving times. An optimal solution is depicted in Figure \ref{fig:outlier}, where the truck route is the sequence $(0,10,9,8,1,5,6,7,4,11)$, while the two sorties are $\langle 0,2,1 \rangle$ and $\langle 1,3,11 \rangle$. This instance is hard to solve due to the large number of possible feasible sorties and to the closeness of customers  among them and to the depot. This allows to have many possible feasible solutions with similar cost, if not with the very same one, e.g. (0,4,,6,5,1,8,9,10,11),$\langle 0,3,1 \rangle$, and $\langle 1,2,11 \rangle$. We believe that this impacts on the convergence speed of the algorithm.}
\end{sloppypar}

\subsection{Medium Instances}\label{mm}
Since all small size instances could be solved to optimality, we challenged our methods by testing them with instances having more than 10 customers. In Murray and Chu \cite{murray2015flying}, the authors propose a set of 120 instances with 20 customers for the PDSTSP that could be easily adapted to the FSTSP. We used the same values as for the small size instances for $\sigma^L = \sigma^R = 1$ and we set the endurance to 20 and 40 unit times. A total of 120 instances with $E=20$ and 120 instances with $E=40$ were thus obtained.

In Table \ref{tab:m20} and Table \ref{tab:m40} we report the results achieved by the most promising models (2IF and 2IF-BC) and by RRLS on the set of medium instances, with endurance $E=20$ and $E=40$, respectively. For 2IF and 2IF-BC we report the upper bound retrieved with a maximum computation time of 3600 seconds, and the optimality gap (calculated as $opt gap\% = 100 \cdot (UB - LB)/LB$, where UB and LB are the upper and lower bounds found by the method) at the end of the given time. For the heuristic RRLS, we report the gap with respect to the best upper bound ($Best_{UB}$) retrieved by either 2IF or 2IF-BC and calculated as $gap UB\% = 100 \cdot (UB - Best_{UB})/Best_{UB}$. A maximum computation time of 720 seconds was considered for RRLS in this case.
{At the end of the tables we report the average of $opt gap\%$ and $gap UB\%$, the average computing times, and the number of optima obtained by the exact methods, while for the RRLS we report the number of solutions with the same value of the best of the upper bounds and the number of improved solutions.}

Table \ref{tab:m20} suggests that 2IF and 2IF-BC are able to obtain good results on the instances of this set when $E=20$. In particular, 2IF-BC performs on average better both in terms of best heuristic solution retrieved and optimality gap, although in some cases 2IF finds better heuristic solutions. Model 2IF-BC is able to close 50 instances, while  2IF stops at 28. RRLS performs well on all the instances in the given time, with small gaps and 5 heuristic solutions improving those of the other methods. The results of Table \ref{tab:m40}, where the endurance $E$ is increased to 40, show a different behaviour. The exact methods have much larger optimality gaps, and {only 2IF-BC is able to solve merely 2 instances, while 2IF fails in all cases}. Comparing the two approaches, now for several instances 2IF has better heuristic results and gaps than 2IF-BC, and an explanation for such a behaviour is the same already provided in Section \ref{sec:results_small}: by increasing the freedom of the drone, 2IF-BC requires to separate much more constraints. {On the other hand, 2IF-BC always provides better lower bounds than 2IF.} When analyzing the results provided by RRLS, it can now be observed that there is an improvement over those retrieved by the other methods in 89 cases, and the improvement is often substantial. This can be interpreted as a signal that when the exact methods reach their limit, RRLS can represent a valid method to retrieve good quality solutions.

\subsection{Large Instances}
Medium instances appear to be already challenging for the proposed exact methods, however we intended to test the algorithms on larger instances, especially to understand how RRLS scales up.
To do so we used the benchmark instances based on the TSPLIB that Mbiadou Saleu et al. \cite{mbiadou2018iterative} proposed for the PDSTSP that were easily adapted to the FSTSP. We selected those between 48 and 101 customers, with several percentage of drone eligible nodes (where 0\% means that no node can be visited by drone resulting in a TSP instance, and 100\% means that all nodes can be served by the drone), with several drone speeds, and with different depot locations. A total of 33 instances is considered.

We compare the best exact algorithms, 2IF and 2IF-BC (run for at most 3600 seconds), and the matheuristic RRLS (with a maximum time of 720 seconds). Results are presented in Table \ref{tab:l} where an instance is characterized by the starting TSPLIB problem, the percentage of drone-eligible nodes (a), the drone speed (b) and the depot location (c). {The drone speed (b) gives the ratio between the drone speed and the truck one, while (c) takes value 1 when the depot is located at the center of all the customers, and value 2 when located at the left‐bottom corner of the considered region}. For the algorithms, the same information reported for the tests of Section \ref{mm} are considered.

From Table \ref{tab:l} one can see that the percentage of drone eligible nodes is the feature with the highest impact on the convergence of the exact algorithm, that solves more efficiently instances with a smaller drone eligible percentage. 
In general, the optimality gaps of the exact methods are very high, indicating that their limit has been reached and probably passed. It is interesting to observe how 2IF is sometimes performing better than 2IF-BC: in a few circumstances it is even able to retrieve heuristic solutions while 2IF-BC is not.  The heuristic RRLS is consistently better (apart from outlier cases), normally finding heuristic solutions substantially better than those of the other methods. This suggests that RRLS can be regarded as a robust method, able to cope with challenging instances.

\begin{landscape}
\begin{table}[htbp]
  \centering
  \tiny
  \caption{Results on medium instances from \cite{murray2015flying} with endurance $E=20$.}
    \begin{tabular}{lcrcrrlcrcrrlcrcrr}
    \toprule \\
    \cmidrule(lr){1-6}\cmidrule(lr){7-12}\cmidrule(lr){13-18}
    \multicolumn{1}{c}{Inst.} & \multicolumn{2}{c}{2IF} & \multicolumn{2}{c}{2IF-BC} & \multicolumn{1}{c}{RRLS}& \multicolumn{1}{c}{Inst.} & \multicolumn{2}{c}{2IF} & \multicolumn{2}{c}{2IF-BC} & \multicolumn{1}{c}{RRLS} & \multicolumn{1}{c}{Inst.} & \multicolumn{2}{c}{2IF} & \multicolumn{2}{c}{2IF-BC} & \multicolumn{1}{c}{RRLS}\\
    \multicolumn{1}{c}{} & \multicolumn{2}{c}{(3600s)} & \multicolumn{2}{c}{(3600s)} & \multicolumn{1}{c}{(720s)}& \multicolumn{1}{c}{} & \multicolumn{2}{c}{(3600s)} & \multicolumn{2}{c}{(3600s)} & \multicolumn{1}{c}{(720s)} & \multicolumn{1}{c}{} & \multicolumn{2}{c}{(3600s)} & \multicolumn{2}{c}{(3600s)} & \multicolumn{1}{c}{(720s)}\\
    \cmidrule(lr){1-1} \cmidrule(lr){2-3}\cmidrule(lr){4-5}\cmidrule(lr){6-6}\cmidrule(lr){7-7} \cmidrule(lr){8-9}\cmidrule(lr){10-11}\cmidrule(lr){12-12} \cmidrule(lr){13-13} \cmidrule(lr){14-15}\cmidrule(lr){16-17}\cmidrule(lr){18-18}
    \!\!\!	&UB&opt&UB&opt&gap& \!\!\! 	&UB&opt&UB&opt&gap& \!\!\! 	&UB&opt&UB&opt&gap\\
    \!\!\!   &   & gap\%  &   & gap\% & UB \%&  \!\!\!   &   & gap\%    &   & gap\%     & UB \%&  \!\!\!   &   & gap\%    &   & gap\%     & UB \%\\
                          \cmidrule(lr){1-1} \cmidrule(lr){2-3}\cmidrule(lr){4-5}\cmidrule(lr){6-6}\cmidrule(lr){7-7} \cmidrule(lr){8-9}\cmidrule(lr){10-11}\cmidrule(lr){12-12} \cmidrule(lr){13-13} \cmidrule(lr){14-15}\cmidrule(lr){16-17}\cmidrule(lr){18-18}
4847	&	267.05	&	5.25	&	267.05	&	2.29	&	1.11	&	5025	&	131.43	&	9.89	&	131.43	&	2.58	&	2.77	&	5154	&	124.82	&	18.61	&	123.34	&	13.76	&	0.00	\\
4849	&	248.30	&	0	&	248.30	&	0	&	0.00	&	5027	&	115.31	&	24.23	&	120.27	&	34.89	&	1.14	&	5156	&	124.46	&	9.09	&	124.46	&	0	&	0.00	\\
4853	&	232.87	&	0	&	232.87	&	0	&	0.14	&	5030	&	117.41	&	19.64	&	117.96	&	24.13	&	2.79	&	5159	&	145.92	&	5.21	&	145.79	&	0	&	1.31	\\
4856	&	253.33	&	0	&	253.33	&	0	&	0.00	&	5032	&	117.55	&	32.17	&	117.55	&	18.46	&	2.45	&	5201	&	148.02	&	0	&	148.02	&	0	&	0.00	\\
4858	&	240.63	&	9.35	&	240.63	&	3.22	&	0.31	&	5034	&	105.10	&	25.93	&	105.18	&	16.97	&	0.00	&	5203	&	138.92	&	5.82	&	138.59	&	0	&	0.72	\\
4902	&	242.32	&	0	&	242.32	&	0	&	0.00	&	5036	&	124.33	&	27.79	&	124.33	&	28.77	&	1.07	&	5205	&	134.78	&	17.29	&	135.80	&	4.35	&	3.34	\\
4907	&	239.28	&	0	&	239.28	&	0	&	0.58	&	5039	&	130.91	&	19.20	&	134.23	&	19.87	&	3.50	&	5207	&	121.47	&	8.44	&	121.47	&	0	&	0.00	\\
4909	&	222.88	&	0	&	222.88	&	0	&	1.42	&	5041	&	125.57	&	18.31	&	125.32	&	8.87	&	3.82	&	5209	&	135.92	&	0	&	135.92	&	0	&	2.23	\\
4912	&	267.62	&	0	&	267.62	&	0	&	0.51	&	5044	&	121.87	&	24.38	&	125.40	&	38.31	&	0.97	&	5212	&	137.67	&	0	&	137.67	&	0	&	3.07	\\
4915	&	259.39	&	0	&	259.46	&	0	&	0.02	&	5047	&	112.84	&	7.64	&	113.92	&	8.98	&	0.96	&	5214	&	126.59	&	10.20	&	126.25	&	0	&	5.01	\\
4917	&	173.97	&	23.14	&	173.97	&	21.47	&	0.00	&	5049	&	198.28	&	3.01	&	197.76	&	0	&	0.66	&	5216	&	105.96	&	29.25	&	101.07	&	22.06	&	4.97	\\
4920	&	170.05	&	16.44	&	172.03	&	15.57	&	0.17	&	5051	&	180.62	&	3.82	&	180.62	&	0	&	0.49	&	5218	&	116.00	&	11.13	&	115.82	&	5.79	&	0.16	\\
4922	&	169.60	&	20.56	&	169.60	&	25.18	&	0.00	&	5053	&	176.51	&	0	&	176.51	&	0	&	0.31	&	5220	&	119.04	&	11.91	&	119.04	&	0	&	1.66	\\
4924	&	160.01	&	19.72	&	159.57	&	12.65	&	0.00	&	5055	&	177.29	&	0	&	177.29	&	0	&	0.00	&	5223	&	94.64	&	19.18	&	94.59	&	3.64	&	0.05	\\
4926	&	155.98	&	20.53	&	155.98	&	14.67	&	0.00	&	5057	&	180.77	&	9.12	&	180.77	&	27.43	&	0.99	&	5225	&	133.07	&	14.72	&	129.72	&	6.72	&	2.30	\\
4928	&	167.61	&	16.19	&	166.27	&	2.07	&	0.80	&	5059	&	150.81	&	10.73	&	150.82	&	0	&	3.12	&	5227	&	116.47	&	18.72	&	116.19	&	0	&	2.42	\\
4931	&	172.49	&	10.17	&	172.49	&	4.22	&	2.32	&	5102	&	165.71	&	11.11	&	165.49	&	6.50	&	0.13	&	5229	&	94.26	&	16.56	&	99.14	&	24.19	&	2.50	\\
4933	&	159.39	&	14.39	&	159.39	&	7.09	&	1.74	&	5104	&	181.61	&	11.33	&	181.89	&	6.37	&	1.29	&	5231	&	100.54	&	23.03	&	99.73	&	7.01	&	0.81	\\
4935	&	178.27	&	18.45	&	177.39	&	12.39	&	0.54	&	5106	&	159.49	&	11.94	&	158.49	&	0	&	0.96	&	5233	&	111.62	&	10.26	&	113.51	&	8.97	&	1.68	\\
4937	&	173.55	&	5.00	&	173.55	&	1.05	&	1.08	&	5108	&	172.12	&	0	&	172.12	&	0	&	0.16	&	5235	&	123.55	&	9.95	&	118.89	&	0.29	&	0.00	\\
4939	&	201.03	&	0	&	201.03	&	0	&	0.52	&	5110	&	135.49	&	17.43	&	135.43	&	0	&	3.68	&	5238	&	79.39	&	27.92	&	79.39	&	32.14	&	0.00	\\
4941	&	253.08	&	0	&	253.08	&	0	&	0.42	&	5112	&	133.09	&	17.29	&	131.23	&	0	&	1.17	&	5240	&	88.71	&	35.63	&	89.79	&	25.82	&	-0.57	\\
4944	&	247.03	&	5.82	&	247.50	&	4.88	&	0.40	&	5115	&	127.53	&	13.16	&	127.44	&	6.91	&	0.07	&	5242	&	88.20	&	9.09	&	85.65	&	3.25	&	0.00	\\
4946	&	237.21	&	0	&	237.21	&	0	&	0.23	&	5117	&	135.94	&	30.79	&	130.36	&	10.31	&	1.91	&	5244	&	86.81	&	0	&	86.81	&	0	&	0.00	\\
4948	&	258.06	&	0	&	258.06	&	0	&	0.28	&	5119	&	120.16	&	22.43	&	118.97	&	0	&	5.52	&	5246	&	75.53	&	29.28	&	74.56	&	17.16	&	0.37	\\
4950	&	239.46	&	4.57	&	240.99	&	0	&	0.64	&	5121	&	131.84	&	12.56	&	132.68	&	11.64	&	2.07	&	5248	&	83.10	&	37.30	&	89.03	&	29.37	&	1.60	\\
4952	&	218.09	&	0	&	218.09	&	0	&	0.00	&	5123	&	121.95	&	22.70	&	121.95	&	10.45	&	1.69	&	5250	&	81.93	&	23.71	&	81.93	&	0	&	0.40	\\
4954	&	261.06	&	6.51	&	261.06	&	0	&	0.22	&	5125	&	132.32	&	21.87	&	130.96	&	7.38	&	2.85	&	5252	&	87.38	&	28.29	&	90.37	&	23.76	&	0.27	\\
4957	&	253.67	&	14.55	&	252.28	&	8.09	&	1.17	&	5127	&	132.24	&	7.25	&	132.24	&	0	&	2.24	&	5255	&	82.50	&	21.61	&	82.50	&	7.10	&	3.37	\\
4959	&	249.92	&	0	&	249.92	&	0	&	0.42	&	5130	&	126.49	&	17.79	&	126.50	&	9.37	&	1.87	&	5257	&	78.98	&	30.57	&	79.43	&	39.68	&	1.10	\\
5001	&	115.88	&	7.92	&	115.50	&	0	&	0.87	&	5132	&	106.36	&	31.47	&	106.73	&	28.95	&	1.77	&	5306	&	100.46	&	0	&	100.46	&	0	&	1.42	\\
5003	&	173.54	&	5.29	&	173.54	&	3.47	&	-0.35	&	5134	&	103.57	&	15.05	&	102.57	&	15.01	&	0.00	&	5310	&	92.79	&	23.47	&	92.47	&	7.92	&	0.00	\\
5006	&	155.39	&	0	&	155.39	&	0	&	2.31	&	5136	&	99.02	&	27.29	&	98.91	&	25.67	&	0.00	&	5312	&	83.59	&	0	&	83.60	&	0	&	0.94	\\
5008	&	159.74	&	0	&	159.74	&	0	&	0.47	&	5138	&	91.97	&	27.13	&	92.49	&	16.38	&	0.72	&	5321	&	101.49	&	0	&	101.49	&	0	&	2.01	\\
5010	&	146.48	&	14.10	&	145.48	&	7.84	&	1.48	&	5141	&	95.53	&	29.17	&	96.74	&	31.70	&	0.22	&	5324	&	104.03	&	14.02	&	105.22	&	13.96	&	-2.00	\\
5012	&	172.65	&	10.51	&	172.40	&	3.42	&	2.96	&	5143	&	98.61	&	30.06	&	99.24	&	27.44	&	-0.94	&	5330	&	118.45	&	0.08	&	118.45	&	0	&	0.91	\\
5015	&	172.67	&	7.27	&	172.67	&	4.01	&	1.79	&	5145	&	95.21	&	25.33	&	95.09	&	33.49	&	0.00	&	5334	&	101.53	&	6.57	&	102.02	&	0	&	1.43	\\
5017	&	171.05	&	14.26	&	166.47	&	7.13	&	2.33	&	5148	&	90.58	&	19.15	&	94.30	&	14.93	&	0.99	&	5336	&	104.46	&	0 	&	104.46	&	0	&	0.00	\\
5020	&	156.07	&	19.33	&	155.92	&	0	&	1.33	&	5150	&	82.75	&	20.34	&	85.87	&	28.87	&	-0.21	&	5345	&	114.19	&	0 	&	114.19	&	0	&	0.10	\\
5022	&	146.69	&	12.71	&	146.21	&	3.41	&	1.96	&	5152	&	92.04	&	42.56	&	90.58	&	30.55	&	2.09	&	5351	&	115.88	&	18.63	&	115.99	&	1.88	&	1.91	\\
 \cmidrule(lr){1-1} \cmidrule(lr){2-3}\cmidrule(lr){4-5}\cmidrule(lr){6-6}\cmidrule(lr){7-7} \cmidrule(lr){8-9}\cmidrule(lr){10-11}\cmidrule(lr){12-12} \cmidrule(lr){13-13} \cmidrule(lr){14-15}\cmidrule(lr){16-17}\cmidrule(lr){18-18}
Avg. time	&	&	&	&	&	&	&	&	&	&	&	&	& & 3048.95	& 	& 2465.94	& 172.02	\\
Avg. \%gap	&	&	&	&	&	&	&	&	&	&	&	&	& & 12.96	& 	& 8.38	&  0.71	\\
\#opt (impr)	&	&	&	&	&	&	&	&	&	&	&	& &	& 28 & 	& 50 &  24 (15)	\\
\bottomrule
    \end{tabular}%
  \label{tab:m20}%
\end{table}%
\end{landscape}
\begin{landscape}
\begin{table}[htbp]
  \centering
  \tiny
  \caption{Results on medium  instances from \cite{murray2015flying} with endurance $E=40$.}
    \begin{tabular}{lcrcrrlcrcrrlcrcrr}
    \toprule \\
    \cmidrule(lr){1-6}\cmidrule(lr){7-12}\cmidrule(lr){13-18}
    \multicolumn{1}{c}{Inst.} & \multicolumn{2}{c}{2IF} & \multicolumn{2}{c}{2IF-BC} & \multicolumn{1}{c}{RRLS}& \multicolumn{1}{c}{Inst.} & \multicolumn{2}{c}{2IF} & \multicolumn{2}{c}{2IF-BC} & \multicolumn{1}{c}{RRLS} & \multicolumn{1}{c}{Inst.} & \multicolumn{2}{c}{2IF} & \multicolumn{2}{c}{2IF-BC} & \multicolumn{1}{c}{RRLS}\\
    \multicolumn{1}{c}{} & \multicolumn{2}{c}{(3600s)} & \multicolumn{2}{c}{(3600s)} & \multicolumn{1}{c}{(720s)}& \multicolumn{1}{c}{} & \multicolumn{2}{c}{(3600s)} & \multicolumn{2}{c}{(3600s)} & \multicolumn{1}{c}{(720s)} & \multicolumn{1}{c}{} & \multicolumn{2}{c}{(3600s)} & \multicolumn{2}{c}{(3600s)} & \multicolumn{1}{c}{(720s)}\\
    \cmidrule(lr){1-1} \cmidrule(lr){2-3}\cmidrule(lr){4-5}\cmidrule(lr){6-6}\cmidrule(lr){7-7} \cmidrule(lr){8-9}\cmidrule(lr){10-11}\cmidrule(lr){12-12} \cmidrule(lr){13-13} \cmidrule(lr){14-15}\cmidrule(lr){16-17}\cmidrule(lr){18-18}
    \!\!\!	&UB&opt&UB&opt&gap& \!\!\! 	&UB&opt&UB&opt&gap& \!\!\! 	&UB&opt&UB&opt&gap\\
    \!\!\!   &   & gap\%  &   & gap\% & UB \%&  \!\!\!   &   & gap\%    &   & gap\%     & UB \%&  \!\!\!   &   & gap\%    &   & gap\%     & UB \%\\
                          \cmidrule(lr){1-1} \cmidrule(lr){2-3}\cmidrule(lr){4-5}\cmidrule(lr){6-6}\cmidrule(lr){7-7} \cmidrule(lr){8-9}\cmidrule(lr){10-11}\cmidrule(lr){12-12} \cmidrule(lr){13-13} \cmidrule(lr){14-15}\cmidrule(lr){16-17}\cmidrule(lr){18-18}
4847	&	255.60	&	10.70	&	295.62	&	26.96	&	-11.33	&	5025	&	123.65	&	32.26	&	119.20	&	42.82	&	6.26	&	5154	&	108.99	&	37.06	&	106.64	&	29.87	&	0.56	\\
4849	&	228.29	&	38.59	&	225.15	&	22.37	&	1.45	&	5027	&	114.12	&	29.30	&	112.15	&	33.23	&	1.07	&	5156	&	108.92	&	21.98	&	116.86	&	16.20	&	-6.20	\\
4853	&	218.60	&	18.68	&	236.21	&	29.14	&	-6.16	&	5030	&	102.69	&	24.57	&	108.99	&	39.97	&	-3.95	&	5159	&	128.34	&	36.11	&	125.41	&	17.76	&	-4.29	\\
4856	&	241.18	&	21.12	&	237.03	&	23.90	&	4.52	&	5032	&	107.53	&	38.39	&	120.20	&	74.15	&	-13.43	&	5201	&	140.30	&	23.47	&	144.27	&	20.53	&	-1.60	\\
4858	&	215.63	&	13.81	&	226.81	&	21.56	&	-3.87	&	5034	&	102.63	&	24.87	&	112.23	&	46.66	&	-8.55	&	5203	&	124.20	&	19.56	&	124.54	&	16.32	&	0.96	\\
4902	&	226.62	&	27.12	&	234.69	&	26.28	&	-3.29	&	5036	&	112.34	&	20.40	&	117.96	&	29.32	&	-4.86	&	5205	&	119.79	&	40.06	&	128.17	&	57.45	&	-5.55	\\
4907	&	201.19	&	24.85	&	218.93	&	25.60	&	-8.38	&	5039	&	126.96	&	44.84	&	127.99	&	73.84	&	-5.56	&	5207	&	113.90	&	31.05	&	115.82	&	32.57	&	-1.66	\\
4909	&	216.04	&	25.93	&	217.36	&	15.20	&	-0.60	&	5041	&	114.55	&	17.12	&	115.54	&	17.75	&	0.42	&	5209	&	123.69	&	49.66	&	139.95	&	78.16	&	-10.98	\\
4912	&	238.62	&	13.41	&	238.04	&	1.83	&	4.33	&	5044	&	115.94	&	23.26	&	123.43	&	40.25	&	-4.41	&	5212	&	136.40	&	30.47	&	138.88	&	25.61	&	-1.78	\\
4915	&	241.13	&	21.26	&	234.42	&	13.26	&	2.03	&	5047	&	106.11	&	20.30	&	107.45	&	25.77	&	-1.25	&	5214	&	123.05	&	47.66	&	135.18	&	42.11	&	-6.58	\\
4917	&	165.31	&	20.89	&	165.62	&	11.23	&	-0.19	&	5049	&	185.73	&	22.13	&	183.51	&	20.39	&	4.86	&	5216	&	95.60	&	19.68	&	104.53	&	33.47	&	-11.04	\\
4920	&	162.89	&	38.28	&	171.18	&	43.17	&	-7.40	&	5051	&	165.69	&	35.87	&	165.99	&	15.32	&	1.47	&	5218	&	97.61	&	21.93	&	120.89	&	57.59	&	-14.39	\\
4922	&	168.76	&	31.96	&	172.41	&	41.90	&	-2.75	&	5053	&	146.59	&	34.96	&	145.54	&	26.26	&	7.68	&	5220	&	119.27	&	16.75	&	131.26	&	28.62	&	-8.86	\\
4924	&	158.98	&	30.11	&	158.98	&	19.43	&	0.00	&	5055	&	173.69	&	13.72	&	174.09	&	10.36	&	-1.38	&	5223	&	93.64	&	22.49	&	107.83	&	48.41	&	-14.52	\\
4926	&	152.24	&	45.63	&	165.67	&	75.59	&	-8.31	&	5057	&	174.03	&	38.98	&	173.67	&	16.41	&	0.29	&	5225	&	126.57	&	10.02	&	130.71	&	14.86	&	-2.28	\\
4928	&	153.88	&	40.19	&	185.82	&	83.65	&	-17.19	&	5059	&	134.83	&	16.19	&	137.52	&	14.06	&	-1.96	&	5227	&	105.35	&	18.49	&	126.96	&	47.85	&	-16.94	\\
4931	&	167.35	&	31.22	&	169.29	&	12.96	&	-1.14	&	5102	&	164.01	&	28.11	&	168.79	&	20.93	&	-2.77	&	5229	&	94.26	&	17.15	&	103.01	&	31.69	&	-6.13	\\
4933	&	157.19	&	37.26	&	171.11	&	63.46	&	-8.54	&	5104	&	180.25	&	27.28	&	178.67	&	17.80	&	3.20	&	5231	&	99.17	&	29.33	&	107.87	&	45.70	&	-7.91	\\
4935	&	161.97	&	40.02	&	164.11	&	60.83	&	-2.09	&	5106	&	144.37	&	33.42	&	158.17	&	48.77	&	-8.73	&	5233	&	108.03	&	16.81	&	114.23	&	27.10	&	-0.65	\\
4937	&	161.10	&	35.55	&	159.44	&	9.91	&	3.10	&	5108	&	165.22	&	39.73	&	162.21	&	10.00	&	3.21	&	5235	&	98.96	&	19.80	&	106.96	&	32.79	&	0.77	\\
4939	&	181.19	&	43.31	&	184.62	&	14.17	&	-2.38	&	5110	&	129.07	&	25.14	&	129.37	&	29.49	&	3.38	&	5238	&	79.09	&	28.75	&	81.85	&	40.17	&	-4.16	\\
4941	&	243.11	&	6.90	&	241.90	&	4.71	&	0.91	&	5112	&	130.16	&	25.21	&	141.98	&	33.42	&	-12.43	&	5240	&	83.27	&	31.17	&	88.80	&	49.17	&	-5.59	\\
4944	&	237.28	&	23.34	&	236.00	&	9.59	&	-0.33	&	5115	&	126.90	&	21.09	&	129.62	&	27.22	&	-3.48	&	5242	&	85.65	&	6.77	&	89.19	&	11.69	&	-3.97	\\
4946	&	225.93	&	26.94	&	220.39	&	0.00	&	2.10	&	5117	&	124.10	&	27.64	&	141.94	&	55.23	&	-10.63	&	5244	&	85.55	&	36.27	&	86.81	&	55.29	&	-1.72	\\
4948	&	244.13	&	25.27	&	246.66	&	24.25	&	2.26	&	5119	&	114.18	&	28.82	&	119.75	&	44.55	&	-3.40	&	5246	&	63.74	&	18.28	&	64.57	&	20.56	&	-1.78	\\
4950	&	232.51	&	19.73	&	229.48	&	14.20	&	-0.65	&	5121	&	118.22	&	7.72	&	122.93	&	12.67	&	-2.84	&	5248	&	83.97	&	39.64	&	92.00	&	76.40	&	-12.53	\\
4952	&	213.30	&	23.60	&	211.67	&	8.56	&	1.33	&	5123	&	115.40	&	19.98	&	136.71	&	43.60	&	-15.41	&	5250	&	81.73	&	24.78	&	85.97	&	36.69	&	-2.99	\\
4954	&	246.41	&	21.35	&	240.59	&	14.62	&	2.84	&	5125	&	120.17	&	25.42	&	130.54	&	36.29	&	-3.11	&	5252	&	85.33	&	21.37	&	97.95	&	51.38	&	-10.56	\\
4957	&	231.61	&	32.32	&	231.61	&	26.76	&	2.57	&	5127	&	126.66	&	13.07	&	140.39	&	26.20	&	-7.68	&	5255	&	73.92	&	27.81	&	75.07	&	39.46	&	-1.57	\\
4959	&	239.87	&	30.49	&	231.95	&	0.00	&	3.89	&	5130	&	120.53	&	25.07	&	139.90	&	53.94	&	-13.85	&	5257	&	79.75	&	28.75	&	90.93	&	61.10	&	-13.62	\\
5001	&	114.09	&	4.37	&	120.25	&	13.18	&	-5.12	&	5132	&	107.47	&	31.23	&	109.31	&	46.82	&	-2.68	&	5306	&	72.37	&	17.01	&	73.06	&	21.65	&	-0.95	\\
5003	&	168.80	&	6.01	&	166.92	&	5.21	&	-1.21	&	5134	&	102.13	&	14.98	&	102.57	&	17.10	&	-1.30	&	5310	&	92.79	&	41.22	&	99.24	&	75.28	&	-8.42	\\
5006	&	135.92	&	15.77	&	135.92	&	12.90	&	1.86	&	5136	&	93.99	&	19.80	&	94.12	&	21.57	&	-0.28	&	5312	&	82.20	&	39.29	&	108.25	&	93.49	&	-24.07	\\
5008	&	155.42	&	38.75	&	156.71	&	19.72	&	-1.18	&	5138	&	87.76	&	21.17	&	102.06	&	44.67	&	-14.00	&	5321	&	90.20	&	26.98	&	79.60	&	6.69	&	1.60	\\
5010	&	135.72	&	19.88	&	150.36	&	35.64	&	-9.19	&	5141	&	96.13	&	30.62	&	97.72	&	44.47	&	-4.06	&	5324	&	91.68	&	41.47	&	92.99	&	54.83	&	0.75	\\
5012	&	161.61	&	22.70	&	160.20	&	23.43	&	-1.51	&	5143	&	97.61	&	29.76	&	98.31	&	33.58	&	-4.52	&	5330	&	96.35	&	41.48	&	101.31	&	31.96	&	-3.14	\\
5015	&	159.02	&	27.80	&	168.58	&	45.42	&	-5.51	&	5145	&	88.13	&	20.32	&	89.01	&	22.33	&	-0.99	&	5334	&	89.68	&	11.70	&	91.47	&	14.35	&	-2.30	\\
5017	&	157.33	&	18.77	&	157.53	&	18.69	&	2.96	&	5148	&	91.82	&	16.43	&	92.94	&	24.83	&	-1.50	&	5336	&	87.41	&	43.82	&	92.22	&	52.71	&	-3.45	\\
5020	&	146.01	&	27.97	&	156.34	&	43.63	&	-6.61	&	5150	&	77.99	&	13.79	&	79.04	&	17.84	&	1.19	&	5345	&	102.00	&	21.97	&	101.45	&	23.94	&	-0.36	\\
5022	&	133.75	&	40.56	&	148.67	&	55.22	&	-10.47	&	5152	&	82.49	&	37.70	&	80.82	&	50.27	&	-0.43	&	5351	&	96.86	&	44.29	&	106.67	&	81.73	&	-8.74	\\
 \cmidrule(lr){1-1} \cmidrule(lr){2-3}\cmidrule(lr){4-5}\cmidrule(lr){6-6}\cmidrule(lr){7-7} \cmidrule(lr){8-9}\cmidrule(lr){10-11}\cmidrule(lr){12-12} \cmidrule(lr){13-13} \cmidrule(lr){14-15}\cmidrule(lr){16-17}\cmidrule(lr){18-18}
Avg. time	&	&	&	&	&	&	&	&	&	&	&	&	& & 3600.00	& 	& 3569.30	& 234.78	\\
Avg. \%gap	&	&	&	&	&	&	&	&	&	&	&	&	& & 26.36	& 	& 32.49	&  -3.69	\\
\#opt (impr)	&	&	&	&	&	&	&	&	&	&	&	& &	& 0 & 	& 2 &  1 (89)	\\
\bottomrule
    \end{tabular}%
  \label{tab:m40}%
\end{table}%
\end{landscape}

\begin{table}[htbp]
  \centering
  \tiny
  \caption{Comparison on large  instances from \cite{mbiadou2018iterative}. (a) \% of drone eligible nodes, (b) drone speed with respect to the truck one, (c) depot location.}
    \begin{tabular}{lrrrrrrrr}
    \toprule
    \multicolumn{4}{c}{Instance} & \multicolumn{2}{c}{2IF (3600s)} & \multicolumn{2}{c}{2IF-BC (3600s)} & \multicolumn{1}{c}{RRLS (720s)} \\
                    \cmidrule(lr){1-4} \cmidrule(lr){5-6}\cmidrule(lr){7-8}\cmidrule(lr){9-9}
    name & (a)& (b)&(c)&UB&opt gap\%&UB& opt gap\%&gap UB\% \\
                          \cmidrule(lr){1-4} \cmidrule(lr){5-6}\cmidrule(lr){7-8}\cmidrule(lr){9-9}
      att48 & 20 & 2 & 1 & 42150.00 & 16.80 & \multicolumn{1}{r}{45606.00} & \multicolumn{1}{r}{26.28} & -5.65 \\
    att48 & 40 & 2 & 1 & 49142.00 & 66.14 & \multicolumn{1}{r}{47806.00} & \multicolumn{1}{r}{58.94} & -27.41 \\
    att48 & 60 & 2 & 1 & 72060.00 & 164.73 & \multicolumn{1}{r}{54194.00} & \multicolumn{1}{r}{96.20} & -36.08 \\
    att48 & 80 & 1 & 1 & 53934.00 & 158.32 & \multicolumn{1}{r}{61632.58} & \multicolumn{1}{r}{193.76} & -36.57 \\
    att48 & 80 & 2 & 1 & 82572.00 & 293.14 & \multicolumn{1}{r}{55214.00} & \multicolumn{1}{r}{165.08} & -39.56 \\
    att48 & 80 & 2 & 2 & 65064.00 & 293.23 & \multicolumn{1}{r}{98426.00} & \multicolumn{1}{r}{502.13} & -48.65 \\
    att48 & 80 & 3 & 1 & 75612.00 & 258.99 & \multicolumn{1}{r}{74070.00} & \multicolumn{1}{r}{253.60} & -54.28 \\
    att48 & 80 & 4 & 1 & 53982.00 & 160.11 & \multicolumn{1}{r}{50206.00} & \multicolumn{1}{r}{138.49} & -30.98 \\
    att48 & 80 & 5 & 1 & 47488.00 & 127.72 & \multicolumn{1}{r}{74888.00} & \multicolumn{1}{r}{256.95} & -26.57 \\
    att48 & 100 & 2 & 1 & 100079.33 & 635.78 & \multicolumn{1}{r}{76460.00} & \multicolumn{1}{r}{473.25} & -54.31 \\
    berlin52 & 20 & 2 & 1 & 9350.00 & 0.00 & \multicolumn{1}{r}{9350.00} & \multicolumn{1}{r}{0.00} & 0.27 \\
    berlin52 & 40 & 2 & 1 & 11865.00 & 52.78 & \multicolumn{1}{r}{8386.15} & \multicolumn{1}{r}{2.77} & 3.79 \\
    berlin52 & 60 & 2 & 1 & 17980.00 & 172.49 & \multicolumn{1}{r}{10770.00} & \multicolumn{1}{r}{59.86} & -26.06 \\
    berlin52 & 80 & 1 & 1 & 24121.04 & 515.26 & \multicolumn{1}{r}{11940.00} & \multicolumn{1}{r}{205.76} & -36.02 \\
    berlin52 & 80 & 2 & 1 & 22162.35 & 468.27 & \multicolumn{1}{r}{10420.00} & \multicolumn{1}{r}{162.80} & -28.85 \\
    berlin52 & 80 & 2 & 2 & 10850.00 & 123.65 & \multicolumn{1}{r}{14930.00} & \multicolumn{1}{r}{207.20} & -28.52 \\
    berlin52 & 80 & 3 & 1 & 15224.81 & 266.63 & \multicolumn{1}{r}{10180.00} & \multicolumn{1}{r}{158.70} & -27.52 \\
    berlin52 & 80 & 4 & 1 & 20937.24 & 438.23 & \multicolumn{1}{r}{11405.00} & \multicolumn{1}{r}{168.65} & -35.73 \\
    berlin52 & 80 & 5 & 1 & 20635.00 & 397.23 & \multicolumn{1}{r}{6805.00} & \multicolumn{1}{r}{69.63} & 8.38 \\
    berlin52 & 100 & 2 & 1 & 17042.89 & 735.95 & \multicolumn{1}{r}{15003.79} & \multicolumn{1}{r}{614.99} & -47.80 \\
    eil101 & 20 & 2 & 1 & 3744.00 & 429.56 & $\infty$ & - & -79.14 \\
    eil101 & 40 & 2 & 1 & 3253.62 & 430.77 & $\infty$ & - & -77.35 \\
    eil101 & 60 & 2 & 1 & 2641.37 & 494.23 & $\infty$ & - & -72.17 \\
    eil101 & 80 & 1 & 1 & 2706.44 & 693.68 & $\infty$ & - & -72.33 \\
    eil101 & 80 & 2 & 1 & 2399.20 & 605.65 & $\infty$ & - & -68.70 \\
    eil101 & 80 & 2 & 2 & 2520.63 & 639.19 & $\infty$ & - & -70.09 \\
    eil101 & 80 & 3 & 1 & 2353.29 & 592.14 & $\infty$ & - & -68.94 \\
    eil101 & 80 & 4 & 1 & 2349.34 & 590.98 & $\infty$ & - & -68.54 \\
    eil101 & 80 & 5 & 1 & 2344.34 & 585.48 & $\infty$ & - & -68.48 \\
    eil101 & 100 & 2 & 1 & 2087.41 & 690.68 & $\infty$ & - & -64.60 \\
        \bottomrule
    \end{tabular}%
  \label{tab:l}%
\end{table}%

\section{Conclusions}\label{sec:conclusions}
In this paper we proposed novel formulations for the Flying Sidekick Traveling Salesman Problem. We also proposed a matheuristic method based on the previous model, and designed for large and difficult instances. Extensive computational tests showed that the best formulation could solve to optimality all the benchmark instances from the literature, with 10 customers, in reasonable time. Several instances with 20 customers could be solved to optimality as well, indicating that the proposed models could be used effectively for instances of small/medium size. For what concerns the matheuristic method we introduced, it has provided robust high quality results for all the benchmarks in relatively short computation times, with remarkable results on the larger instances considered, on which it was able to improve the heuristic solutions provided by the other methods by a large amount.
Future works should explore problems with multiple trucks and drones and taking into account harder constraints.

\bibliographystyle{model1-num-names}
\bibliography{references.bib}






% \end{thebibliography}

\end{document}